\documentclass[11pt]{article}
\usepackage[centertags,intlimits]{amsmath}                     
\usepackage{amsfonts}                     
\usepackage{amsthm}
\usepackage{graphics}
\allowdisplaybreaks[2]
\numberwithin{equation}{section}
\newtheorem{theorem}{Theorem}[section]
\newtheorem{lemma}{Lemma}[section]

\theoremstyle{remark}

\providecommand{\abs}[1]{\lvert #1\rvert}

\newcommand{\nc}{\newcommand}
\nc{\vb}{\mathbf{v}}
\nc{\bx}{\mathbf{x}}
\nc{\by}{\mathbf{y}}
\nc{\bz}{\mathbf{z}}
\nc{\bu}{\mathbf{u}}
\nc{\bv}{\mathbf{v}}
\nc{\ba}{\mathbf{a}}
\nc{\bs}{\mathbf{s}}
\nc{\bq}{\mathbf{q}}
\nc{\bd}{\mathbf{d}}
\nc{\bb}{\mathbf{b}}
\nc{\bc}{\mathbf{c}}
\nc{\bi}{\mathbf{i}}
\nc{\be}{\mathbf{e}}
\nc{\bfr}{\mathbf{r}}
\nc{\bA}{\mathbf{A}}
\nc{\R}{\mathbb R}
\nc{\N}{\mathbb N}
\nc{\C}{\mathbb C}
\nc{\D}{\mathbb D}
\nc{\Z}{\mathbb Z}
\nc{\F}{\mathbf F}
\nc{\bbS}{\mathbb S}
\nc{\B}{\cal B}
\nc{\br}{\bigr}
\nc{\bl}{\bigl}
\nc{\Bl}{\Bigl}
\nc{\Br}{\Bigr}
\nc{\ind}{\mathbf{1}}

\title{A large deviation principle  for join the shortest queue}
\author{Anatolii A. Puhalskii\\
University of Colorado at Denver, Denver, U.S.A. and\\
 Institute of Information
Transmission Problems, Moscow, Russia  \\\\
Alexander A. Vladimirov\\
 Institute of Information
Transmission Problems, Moscow, Russia}
\begin{document}
\maketitle
\sloppy
\begin{abstract}
We consider a join-the-shortest-queue model which is as follows. There
are $K$ single FIFO servers and $M$ arrival processes. The customers from
a given arrival process can be served only by servers from a certain
subset of all servers. The actual destination is the server with the
smallest weighted queue length. The arrival processes are assumed to obey
a large deviation principle while the service is  exponential.
A large deviation principle is established for the queue-length
process.  The action functional is expressed in
terms of solutions to mathematical programming problems.
 The large deviation limit point is
identified as a weak solution to a system of idempotent
equations. Uniqueness of the weak solution is proved by establishing
trajectorial uniqueness.
\end{abstract}
{\renewcommand{\thefootnote}{}
\section{Introduction and summary}
\label{sec:introduction}

\paragraph{Motivation}
\label{sec:motivation}

Queueing models with  
   the join-the-shortest-queue (abbreviated further as JSQ) mechanism
  for routing arriving customers are of interest in various
   application areas, see, e.g.,  Fleming and Simon \cite{FS99}, 
Turner \cite{MR2001i:60165}. In a generic set-up, one considers a
  queueing system consisting   of service  stations arranged in parallel.
The customers arrive exogenously and  join  the
  station that has the least number of customers to serve. If there
  are several  stations that have fewest  customers, 
  then some rule is applied to direct an arriving customer
  to one of these stations. The model
  admits many versions: along with the stream of the
  ``discretionary'' customers that choose the station with the least
  number of customers there may be ``dedicated'' customers who can only
  be served at  specific stations, there may be several streams of 
discretionary customers each of which can join only the stations from a
  certain subset of the set of all stations, 
  the number of servers at the
 service stations can vary from one to infinity, the stations may be
  assigned weights so that the  discretionary
  customers  join the station with the least weighted number of
  customers, upon service
  completion the customers may either leave the system 
 or be routed back, and so on. However,
  all these versions share the common feature that  the arrival
  rates   at the stations depend discontinuously on
   the numbers of customers present. More
  precisely,   the
  rate at which a given station receives customers experiences a jump
  when  the set of the stations with fewest customers changes. 

Thus, JSQ models fall in the category of 
stochastic systems  with discontinuous
  dynamics (or discontinuous stochastic dynamical systems). The 
dynamical properties of such a system
  change abruptly when it enters   a certain domain of the
  state space.  In this paper 
we concern ourselves with a study of the large deviation
  principle (LDP)   for the queue lengths in a JSQ model.
Large deviations of discontinuous stochastic dynamical systems 
have received considerable attention
in the literature,
see  Alanyali and Hajek \cite{AH},
Atar and Dupuis  \cite{MR2001a:60030}, 
Blinovski and Dobrushin \cite{BliDob94}, 
Borovkov and Mogulskii \cite{BorMog01},
Bou\'e, Dupuis and Ellis \cite{MR1736592}, 
Dupuis and Ellis \cite{DupEll92,DupEll97,dupel}, Dupuis, Ellis and
Weiss \cite{DupEllWei91}, Ignatiouk \cite{MR2002i:60054,Ign05}, 
Korostelev and Leonov \cite{KorLeo92a,KorLeo93},
Majewski \cite{Maj04},
 Ramanan and Dupuis \cite{RamDup98}, Shwartz and Weiss
 \cite{ShwWei95}. 
However, available results on the LDP for JSQ 
   are confined to the setting of two stations,
see   Ridder and Shwartz \cite{RidShw05} and references therein. 
Even for that case 
we have not been able to find in the literature an explicit
 formulation of an LDP for the queue length process
 with an indication of the associated action
 functional.  Paths to overflow have been
 studied in  Foley and
 McDonald \cite{MR2002h:60207},  
 Ridder and Shwartz \cite{RidShw05}, and Turner \cite{MR2001i:60165}. 
A  version of the model where the stations are represented by infinite servers
 has been analysed by
 Alanyali and Hajek \cite{AH1} who for the case of two stations
establish an LDP for the queue
 length trajectories and study overflow paths, 
and Turner \cite{MR2001i:60165} who contrasts the results for the
 single-server and infinite-server cases.

\paragraph{Model}
\label{sec:model}
We consider a JSQ model with an arbitrary number of arrival processes
and an arbitrary number of service stations
which are represented by  single servers with the FIFO service
discipline.
 The customers from
a given arrival process can be served only by the servers from a certain
subset of all servers. The actual destination is the server with the
smallest weighted queue length.
The  model thus incorporates both  dedicated and discretionary
flows of customers. 
The arrival   processes are general and are only assumed to
 obey an LDP. The service times  are associated with
 the servers and  are exponentially distributied. (In fact, our main
 result  concerns   a more general setting of autonomous service.) 
Upon service completion, the customers depart from the system.

\paragraph{Methods}
\label{sec:methods}
We adopt the strategy that has proved to be useful for establishing
LDPs for continuous dynamical systems.
It is based on the characterisation of large deviation
relatively compact sequences as exponentially tight ones, see, e.g.,
Puhalskii \cite{Puh91}, \cite{Puh01}. The large deviation limit point
is identified as a weak solution to  idempotent equations. 
The latter equations are obtained as  large deviation limits of the
stochastic equations governing the original dynamical system.

However, this general method has to be modified 
for stochastic systems with discontinuous dynamics, in particular, for
the JSQ model we are concerned with here. 
Building on   the approach of
  Puhalskii \cite{Puh05}, we first 
replace  the original discontinuous  equations
  with certain continuous ones, to which  the limit procedure is applied.
 We  prove trajectorial uniqueness for the
resulting system of idempotent equations. This implies uniqueness of
the large deviation limit point, hence, an LDP. 
The trajectorial uniqueness  is proved by showing the existence of 
 a Lyapunov function.

\paragraph{Contribution}
\label{sec:overview-results}

We establish an LDP  for the queue length process
considered as a random element of the associated Skorohod
 space.  The action functional is of integral form and 
is expressed in terms of solutions to
 mathematical programming problems. 
We also provide some insight into the relation
 between weak and
 trajectorial uniqueness for idempotent equations which is
 instrumental in the proof. Besides, uniqueness of solutions 
for a fluid version of the JSQ model is established.

In  broader perspective, the approaches of this paper and those of
Puhalskii \cite{Puh05} provide new tools 
for the study of large
deviations of stochastic systems with discontinuous dynamics. The
system of idempotent equations that needs to have a weakly unique
solution in order for an LDP to hold has the form of  a system of
differential equations. The latter system  can be viewed as describing a fluid
version of the original queueing model where the arrival, service and
queue length processes are represented by absolutely continuous
functions. One of the implications of the results of this paper is that if, 
given trajectories of exogenous arrival and service processes (``the
inputs''),  the
system of differential equations has a unique solution for the queue
length trajectory (``the output'') 
then weak uniqueness for the associated idempotent
equations holds, so an LDP holds. In short, trajectorial 
uniqueness for the fluid
model implies an LDP. In fact, one can allow the fluid model to have a
unique solution for a certain subset of the set of all ``inputs'',
cf., Puhalskii \cite{Puh05}.
As the application in this paper shows, the method does not restrict
 the number of 
domains of constant dynamics sharing a common boundary,
 while  such  constraints are inherent in
 the existing techniques.

\paragraph{Organisation}
\label{sec:organisation}

The paper is organised as follows. 
The main result  is stated in Section~\ref{sec:large-devi-join}.
Section~\ref{sec:proof} is 
concerned with its proof. The appendix reviews  basics of
idempotent probability and large deviation convergence.

\paragraph{Notation, terminology, conventions}
\label{sec:notat-term-conv}

We will say that a function $\mathbf{I}$ from a metric space
$\mathbb{E}$ to $[0,\infty]$ is an action functional if it is lower compact
that is the sets $\{z\in\mathbb{E}:\,\mathbf{I}(z)\le a\}$ are compact
for $a\in\R_+$ and $\inf_{z\in\mathbb{E}}\mathbf{I}(z)=0$.
 A sequence  
$\{ \mathbf{P}_n ,\, n \in\N \}$ 
of probability measures on the Borel $\sigma$-algebra of  
$\mathbb{E}$ (or a sequence of random elements $\{X_n ,\, n \in\N \}$  
with values in $\mathbb{E}$ and distributions $\mathbf{P}_n$) is said to
obey the LDP for scale $n$ with the  
action functional $\mathbf{I}$  if  
$\limsup_{n \to \infty}  
 n^{-1} \log \mathbf{P}_n (F) \leq - \inf_{z \in   F} \mathbf{I}(z)  
$ for each closed subset $F$ of $ \mathbb{E}$ and  
$\liminf_{n \to \infty}  
 n^{-1} \log \mathbf{P}_n (G) \geq - \inf_{z \in G} \mathbf{I}(z)  
$ for each open subset  $G$ of $ \mathbb{E}$.  

The set of natural numbers is denoted as $\N$, the set of real numbers
is denoted as $\R$ and the 
non-negative halfline is denoted as $\R_+$; for
 $\ell\in\N$, $\R^\ell$ and 
$\R_+^\ell$ denote the cartesian products of $\ell$ copies of $\R$ and
$\R_+$, respectively, with product topology.
Elements of $\R^\ell$
 are considered as column-vectors, $1_\ell$ denotes the element of
 $\R^\ell$ with all the entries equal to one,
 superscript $^T$ is used to denote the transpose of a
 matrix. Inequalities involving vectors are understood componentwise.
We denote as  $\D(\R_+,\R^\ell)$
the Skorohod space of
$\R^\ell$-valued right-continuous with left-hand limits  functions on
$\R_+$. It is assumed to be 
 endowed with the Skorohod
 $J_1$-topology and metrised by a complete separable metric,  see
 Ethier and Kurtz \cite{EthKur86}, 
Liptser and Shiryaev \cite{LipShi89}, Jacod and Shiryaev
\cite{jacshir} for the definition and
properties. The elements of
$\D(\R_+,\R^\ell)$ are denoted with lower-case bold-face Roman characters,
e.g., $\bx=(\bx(t),\,t\in\R_+)$. 

The following conventions are assumed: 
sums and infima over the empty set 
are  equal to $0$ and $\infty$, respectively,
$0/0=0$, and $0\cdot (-\infty)=0\cdot \infty=\infty\cdot 0=0$.
For $\alpha\in\R$ and $\beta\in\R$, we 
denote as
 $\lfloor \alpha\rfloor$   the integer part,
 $\alpha\wedge \beta=\min(\alpha,\beta)$,  $\alpha\vee
 \beta=\max(\alpha,\beta)$, and  $\alpha^+=\alpha\vee0$;
  $\ind(\Gamma)$ denotes
 the indicator function of a set $\Gamma$ that is equal to $1$ on
 $\Gamma$ and is equal to $0$ outside of $\Gamma$.
We use dot to denote differentiation with respect to the
 time variable, so that $\dot{\bx}(t)$
denotes the time-derivative of an absolutely continuous
 function $(\bx(t),\,t\in\R_+)$; if
$\bx(t)$ is a vector, then the notation signifies that 
each component is differentiated; $\bx(t-)$ denotes the left-hand
limit at $t$; ``a.e.'' refers to Lebesgue measure
unless specified otherwise. All the relations involving derivatives
are understood to hold a.e.

\section{The LDP for   join the shortest queue}
\label{sec:large-devi-join}

\paragraph{Dynamics}
\label{sec:model-description}

We provide a more specific model description and turn it into
 equations relating stochastic processes.
There are  $K$ queues, indexed $1$ through $K$, 
each with a single  server.
The servers  are fed by  $M$ arrival processes indexed $1$ through
$M$. Customers from the $m$-th arrival process can be  served only by the
servers from a nonempty subset  $S_m$ of all servers. More
 specifically,  an arriving
 customer is routed to the  queue with the smallest weighted queue length:
queue $k\in S_m$ is assigned weight $w_{km}$, where $w_{km}>0$, 
 and the customer in question   joins
the queue with the least value of $x_k/w_{km},\,k=1,2,\ldots,K$, where
$x_k$ stands for the number of customers in queue $k$ at the time when the
arrival occurs. In case there are two or more queues with the least
value of $x_k/w_{km}$, the arrival is routed to one of the competing
queues. The actual rule used for breaking ties is of no
consequence for the results below. 
Note that in this set-up dedicated arrivals
are those with one-element  sets $S_m$.

We assume that the servers perform  service autonomously in the
following sense. Each server   is assigned a point process (i.e.,
a piecewise constant nondecreasing process starting at zero 
 with unit jumps) which specifies service completion times: if at the
 moment preceding a
 jump time  there is a customer present in the queue,
 then this customer leaves the queue at the time of the jump.  
Due to the memoryless property of the exponential distribution,
 servers with exponential service times can be considered as servers
 with  autonomous service.

Our set-up   concerns, in fact, a sequence of queueing models as
described above indexed by $n\in\N$. To distinguish
between the models, superscript $n$
will be used in  the
notation for the entities associated with the $n$-th model.
Let $A^n_m=(A^n_m(t),\,t\in\R_+)$ and 
$B^n_k=(B^n_k(t),\,t\in\R_+)$, where $m=1,2,\ldots,M$ and
$k=1,\ldots,K$, be one-dimensional
point processes. It is assumed that 
$A^n_m(0)=B^n_k(0)=0$.
We let  $A^n_m(t)$ model the cumulative number of exogenous 
 arrivals in the $m$-th arrival process  by time $t$ and let
$B^n_m(t)$ model the number of service completions by server $m$
during time $t$ of uninterrupted work of the server.
 All the processes are  defined on a probability space
$(\Omega,\mathcal{F},\mathbf{P})$ and have trajectories in the
associated Skorohod spaces.    
Let $C_k$ denote the set of $m\in\{1,2,\ldots,M\}$ such that 
$k\in S_m$ and $Q^n_k(t)$ denote  the number of customers in
queue $k$ at time $t$. 
 The
process $Q^n=(Q^n(t),\,t\in\R_+)$, where 
$Q^n(t)=(Q^n_k(t),\,k=1,2,\ldots,K)$, 
satisfies the following balance equations,
where $k=1,2,\ldots, K$   and $t\in\R_+$, 
\begin{multline}
  \label{eq:2}
Q^n_k(t)=Q_k^n(0)+\sum_{m\in C_k}
\Bl[\int_0^t \ind\Bl(\frac{Q^n_k(s-)}{w_{km}}
<\min_{\substack{l\in S_m:\\l\not=k}}\frac{Q^n_{l}(s-)}{w_{lm}}\Br)\,dA^n_m(s)\\+
\int_0^t \ind\Bl(\frac{Q^n_k(s-)}{w_{km}}
=\min_{\substack{l\in S_m:\\l\not=k}}\frac{Q^n_{l}(s-)}{w_{lm}}\Br)
\alpha^n_{km}(s)\,dA^n_m(s)\Br]-
\int_0^t\ind\bl(Q^n_k(s-)>0\br)\,dB^n_k(s)\,.
\end{multline}
In this equation, 
the random variables $\alpha^n_{km}(s)$ account for the rule adopted for
breaking ties between several  queues of minimum weighted
length, so  one of the  random variables
$\alpha^n_{km}(s),\,k=1,2,\ldots,K,$ equals $1$ while the rest equal $0$.
The integrals on the right of \eqref{eq:2} are well defined by being 
 finite sums. Since the  $A^n_m$ and $B^n_k$ have piecewise constant
 trajectories with  finite numbers of jumps on bounded intervals,
 \eqref{eq:2} admits a unique solution for $Q^n_k(t)$, see Chen and
 Mandelbaum \cite{MR1106809} for more extensive results.

\paragraph{Statement of the LDP}
\label{sec:ldp}
Let us be  given a $[0,\infty]$-valued Borel
 function  $\psi(z)$, where $z\in\R_+^{M+K}$.
 Let function $\mathbf{I}^{A,B}:\,\D(\R_+,\R^{M+K})\to\R_+$ be
 defined  by the  equality
 \begin{equation}
   \label{eq:35}
     \mathbf{I}^{A,B}(\bz)=
\int_0^\infty \psi\bl(  \dot{\bz}(t)\br)\,dt,
\end{equation}
if the function $\bz=(\bz(t),\,t\in\R_+)$ is
  absolutely continuous, componentwise nondecreasing, and
 $\bz(0)=0$, and 
$     \mathbf{I}^{A,B}(\bz)=\infty$ otherwise.
It is assumed that $\mathbf{I}^{A,B}$ is an action functional
 on $\D(\R_+,\R^{M+K})$.
It  follows that $\psi$ is an
 action functional itself. 

For $x=(x_1,\ldots,x_K)\in\R_+^K$ 
and $y=(y_1,\ldots,y_K)\in\R^K$, let $N(x,y)$ denote the set 
of
 $a=(a_1,\ldots,a_M)\in\R_+^M$ and $b=(b_1,\ldots,b_K)\in\R_+^K$
 for which there exist  matrices $e=(e_{km})\in \R_+^{K\times M}$
 and  vectors $d=(d_1,\ldots,d_K)\in\R_+^K$  
such that   $y=e1_M-d$,  $e^T 1_K\le a$, $d\le b$,
 $e_{km}=0$ if either $k\not\in S_m$ or
 $x_k/w_{km}>\min_{l\in S_m}x_l/w_{lm}$, and $d_k=a_k$ if $x_k>0$.
Let 
\begin{equation}
  \label{eq:12}
L(x,y)=\inf_{(a,b)\in N(x,y)}
 \psi(a,b).   
\end{equation}
Since $\psi$ is lower compact and $N(x,y)$ is closed, the infimum
above is attained provided $N(x,y)$ is nonempty.
 Besides, an easy compactness argument shows that
the function $L$ is lower semicontinuous. In
particular, it is Borel measurable.

For $q_0\in \R^K$ and $\bq\in \D(\R_+,\R^K)$, we define
\begin{equation*}
  \mathbf{I}^Q_{q_0}(\mathbf{q})=
 \int_0^\infty L(\mathbf{q}(t),\dot{\mathbf{q}}(t))\,dt\,
\end{equation*}
if $q_0\in\R_+^K$, the function $\mathbf{q}$ is componentwise 
nonnegative, absolutely continuous and
$\bq(0)=q_0$, and $\mathbf{I}^Q_{q_0}(\mathbf{q})=\infty$ otherwise.

Let $\overline{Q}^n_k(t)=Q^n_k(nt)/n$, 
$\overline{A}^n_m(t)=A^n_m(nt)/n$, and
$\overline{B}^n_k(t)=B^n_k(nt)/n$, where $m=1,2,\ldots,M$ and $k=1,2,\ldots,K$.
We introduce the processes
$\overline{Q}^n_k=(\overline{Q}^n_k(t),\,t\in\R_+)$,
$\overline{A}^n_m=\bl(\overline{A}^n_m(t),\,t\in\R_+\br)$, 
$\overline{B}^n_m=\bl(\overline{B}^n_m(t),\,t\in\R_+\br)$, 
$\overline{Q}^n=(\overline{Q}^n_k,\,k=1,2,\ldots,K)$,
$\overline{A}^n=\bl(\overline{A}^n_m,\,m=1,2,\ldots,M\br)$, 
and $\overline{B}^n=(\overline{B}^n_k,\,k=1,2,\ldots,K)$.
\begin{theorem}
  \label{the:short}
 Let, as $n\to\infty$, the sequence
$\{(\overline{Q}^n(0),\overline{A}^n,\overline{B}^n),
\,n\in\N\}$ obey the LDP in
$\R^K\times\D(\R_+,\R^M)\times\D(\R_+,\R^K)$  
for scale $n$ with the action
functional $\mathbf{I}^{Q_0,A,B}(q_0,\ba,\bb)
=\mathbf{I}^{Q_0}(q_0)+\mathbf{I}^{A,B}(\ba,\bb)$, where
$\mathbf{I}^{Q_0}$ is an action functional on $\R^K$.
Then the sequence 
$\{\overline{Q}^n,\,n\in\N\}$ obeys the LDP for scale $n$
 in  $\D(\R_+,\R^K)$   with the action
functional $\mathbf{I}^Q$ defined by the equality
 $\mathbf{I}^Q(\bq)=\mathbf{I}^{Q_0}(\bq(0))+\mathbf{I}^{Q}_{\bq(0)}(\bq)$.
\end{theorem}

\paragraph{Comments and corollaries}
\label{sec:comments-corollaries}

  The fact that $\mathbf{I}^Q$ is an action functional under the
  hypotheses is a part of the assertion of the theorem. 
One can  establish this property without relying on the proof of
  Theorem~\ref{the:short}  if in addition to being an action
  functional the function 
$\psi$ is assumed to be 
convex and of superlinear growth at infinity. The latter conditions
  also ensure that  $\mathbf{I}^{A,B}$ is an
  action functional. The argument is similar to the one in Puhalskii
  \cite{Puh05}. 

  Expression \eqref{eq:12} for the local action functional is
  intuitive. If we interpret $a$  as the vector of
  instantaneous exogenous arrival rates,  $b$, as 
the vector of service rates, $d$, as the vector of customer 
departure rates, and $e$, as the matrix of  rates at which  customers 
 from different exogenous  arrival processes arrive at the  servers, then
   \eqref{eq:12} tells us that
  ``the cost'' for the queue-length vector to change at rate $y$
  given it equals $x$ is obtained by minimising the local action
  functional for the exogenous arrivals and service
   subject to certain ``conservation laws'': the rate
  of the queue length change must equal the difference between the 
  arrival and departure rates at the servers, 
the departure rates cannot exceed the
  service rates, the sum over all the servers of the arrival rates  
  due to the customers orginating from a specific exogenous arrival process
  must equal     the arrival rate of  this process. 
 Besides, the ratio $d_k/b_k$, provided $b_k>0$, can be
  interpreted as the fraction of time that server $k$ is busy serving
  customers, while $e_{km}/a_m$ can be interpreted as the fraction of
  stream $m$ customers directed to server $k$.

The definition of the function $L(x,y)$  implies that it is 
piecewise constant  in $x$ in that 
\begin{equation}
  \label{eq:16}
  L(x,y)=\sum_{I,J}  \ind(x\in  F_{IJ})\,\Psi_{IJ}(y)\,.
\end{equation}
The  summation is over all subsets $I$ of $\{1,2,\ldots,K\}$
 (including the empty set) and all sets $J$ 
 of the form $J=J_1\times\ldots\times J_M$, 
where the  $J_m$, for $m=1,2,\ldots,M$, are nonempty subsets of the
 $S_m$ such that either $J_m\subset I$ or $J_m\cap I=\emptyset$, 
  $F_{IJ}$ denotes  the   subset of  $\R^K_+$ of  elements 
$x=(x_1,\ldots,x_K)$ such that $x_k=0$ if $k\in I$, $x_k>0$ if
 $k\not\in I$,
 $x_k/w_{km}=\min_{l\in S_m}
x_l/w_{lm}$  for $m=1,2,\ldots,M$ if $k\in J_m$, and
 $x_k/w_{km}>\min_{l\in S_m}
x_l/w_{lm}$  for $m=1,2,\ldots,M$ if $k\not\in J_m$. In words,
if  $x\in F_{IJ}$, then set $I$ indexes the entries of $x$ that  are
 equal to zero and  $J_m$ indexes the entries that are the
 smallest weighted fluid queue lengths for the $m$-th fluid 
arrival process. The functions 
  $\Psi_{IJ}$ are uniquely specified by \eqref{eq:16}.
It is easy to see that these functions are 
 action functionals on $\R^K$. 
If, moreover,  the function $\psi$ is convex, then the functions
   $\Psi_{IJ}$ are also convex.
Representation \eqref{eq:16} shows that the $F_{IJ}$ are the
 domains of constant queue length dynamics.

The functions $\Psi_{IJ}$
 can be written down  more
  explicitly if more structure is imposed on the function $\psi$.
If we assume that 
$\psi(a,b)=\psi^{A}(a)+\sum_{k=1}^K\psi^B_k(b_k)$, 
where the functions $\psi^B_k:\,\R_+\to  [0,\infty]$ are lower semicontinuous,
 convex, attain zero, and are not equal to zero identically, 
then minimisation with respect to  $b$ in the
  definition of $\Psi_{IJ}$ yields the representation
  \begin{equation*}
    \Psi_{IJ}(y)=\inf_{(a,d)\in H_{J}^{-1}(y)}
\bl(\psi^{A}(a)+\sum_{k\not\in I}\psi^B_k(d_k)
+\sum_{k\in I}\psi^B_k(d_k)\ind(d_k>\mu_k)\br),
  \end{equation*}
where $\mu_k=\sup\{b_k\in\R_+:\, \psi^B_k(b_k)=0\}$ and 
  $H_{J}(a,d)$, for $a\in\R_+^M$ and $d\in\R_+^K$,
 is the set of $\tilde{y}\in\R^K$
 for which there exist 
matrices $e=(e_{km})\in\R_+^{K\times M}$  
such that   $\tilde{y}=e1_M-d$, $e^T1_K\le a$, and $e_{km}=0$ if
$k\not\in J_m$.

As a consequence of this observation and Theorem~\ref{the:short}, 
 we obtain the following
result for the Markovian setting.  For $\alpha\in\R_+$, we denote
 $  \pi(\alpha)=\alpha\log
 \alpha-\alpha+1$.
\begin{theorem}
  \label{the:short_poi}
Let the  $A^n_m$ and $B^n_k$ be independent Poisson processes
with respective rates $\lambda^n_m$ and $\mu^n_k$, which are
also independent of $Q^n(0)$. Let, as $n\to\infty$,
$\lambda^n_m\to\lambda_m$,  $\mu^n_k\to\mu_k$, and
the sequence
$\{\overline{Q}^n(0),\,n\in\N\}$ obey the LDP  for scale $n$ 
in $\R^K$ with an action
functional $\mathbf{I}^{Q_0}$.
Then the sequence 
$\{\overline{Q}^n,\,n\in\N\}$ obeys the LDP for scale $n$
 in  $\D(\R_+,\R^K)$   with the action
functional  $\mathbf{I}^Q(\bq)=\mathbf{I}^{Q_0}(\bq(0))+
\mathbf{I}^Q_{\bq(0)}(\bq)$. The function $L(x,y)$ is of the form
\eqref{eq:16}, where
\begin{equation*}
       \Psi_{IJ}(y)=\inf_{(a,d)\in H_J^{-1}(y)}
\Bl(\sum_{m=1}^M\pi\bl(\frac{a_m}{\lambda_m}\br)\lambda_m
+\sum_{k\not\in I}\pi\bl(\frac{d_k}{\mu_k}\br)\mu_k
+\sum_{k\in I} \pi\bl(\frac{d_k}{\mu_k}\br)\mu_k
\,\ind\bl(\frac{d_k}{\mu_k}>1\br)\Br).
\end{equation*}
\end{theorem}
For the proof of Theorem~\ref{the:short_poi}, note that 
 by Theorem~2.3 in Puhalskii \cite{Puh94b} the LDP for 
$\{(\overline{Q}^n(0),\overline{A}^n,\overline{B}^n),
\,n\in\N\}$ in the hypotheses of Theorem~\ref{the:short} holds with 
$\psi(a,b)=\sum_{m=1}^M\pi(a_m/\lambda_m)\lambda_m+
\sum_{k=1}^K \pi(b_k/\mu_k)\mu_k$.
\section{Proof of Theorem~\ref{the:short}}
\label{sec:proof}

The proof uses the terminology of large deviation convergence
(abbreviated as LD convergence) and
idempotent probability, which is recapitulated in the appendix, for
more detail see Puhalskii \cite{Puh01}. We start with a proof outline,
where we also  define some concepts which illuminate the connection of LD
convergence and weak convergence.

\paragraph{Weak uniqueness and trajectorial uniqueness}
\label{sec:weak-uniq-traj}
 Let 
$X^n=(\overline{Q}^n(0),
\overline{A}^n,\overline{B}^n)$ and $Y^n=\overline{Q}^n$.
 These random variables assume values 
in the respective metric spaces 
$\mathbb{E}^X=\R^K\times\D(\R_+,\R^M)\times\D(\R_+,\R^K)$ and
$\mathbb{E}^Y=\D(\R_+,\R^K)$. The sequence $X^n$ LD converges in
distribution at
rate $n$ by hypotheses.
We seek to prove that the $Y^n$ LD converge in distribution and find the limit.

Equation \eqref{eq:2}
is transformed   into an equation 
$F_t(X^n,Y^n)=0,\,t\in\R_+$, where the functions $F_t$ are discontinuous.
 This equation enables us to establish that the
sequence $Y^n$ is exponentially tight. Thus, there exists a
subsequence $(X^{n'},Y^{n'})$ that LD converges in distribution to a certain
idempotent process $(X,Y)$. By hypotheses, the idempotent process
 $X$ has idempotent distribution 
$\mathbf{\Pi}^X(\bx)=\exp(-\mathbf{I}^{Q_0,A,B}(\bx))$, where
$\bx\in\mathbb{E}^X$.   
In order to identify the idempotent distribution of $Y$ we would like to
relate $X$ and $Y$ by  equations obtained 
as certain LD limits. Since the original 
 equation $F_t(X^n,Y^n)=0,\,t\in\R_+$
 involves discontinuities, as a preliminary step we replace
 it with an equation $\hat{F}_t(X^n,Y^n,Z^n)=0,\,t\in\R_+$, 
where the $\hat{F}_t$ are
 continuous functions and $Z^n$ are additional random variables, which we call
 ``latent''. They assume values in a metric space $\mathbb{E}^Z$. 
The sequence $(X^n,Y^n,Z^n)$ is still provably exponentially
 tight. Taking an LD limit along a subsequence,  we have
 that if  idempotent variables $(X,Y,Z)$ defined on an idempotent
 probability space $(\Upsilon,\mathbf{\Pi})$
constitute  an LD accumulation point of the
 $(X^n,Y^n,Z^n)$ for LD convergence in distribution, then
 $\hat{F}_t(X,Y,Z)=0$ $\mathbf{\Pi}$-a.e. (Note that one can always take
 $\mathbb{E}^X\times\mathbb{E}^Y\times\mathbb{E}^Z$ as
 $\Upsilon$.)  The idempotent distribution
 of $(X,Y,Z)$ is proved to be
 concentrated on a set  $\mathbb{E}_0\subset
  \mathbb{E}^X\times\mathbb{E}^Y\times\mathbb{E}^Z$ in the sense  that 
$\mathbf{\Pi}((X,Y,Z)\in
\mathbb{E}^X\times\mathbb{E}^Y\times\mathbb{E}^Z\setminus
\mathbb{E}_0)=0$.  
We refer to the idempotent
 distribution of $Y$, which is defined by
 $\mathbf{\Pi}^Y(\by)=\mathbf{\Pi}(Y=\by)$, 
 as a weak solution of the equation $\hat{F}_t(X,Y,Z)=0$. If this
  weak solution is unique, then the idempotent law of $Y$ is specified
  uniquely, so the $Y^n$ LD converge in distribution to $Y$. 
As in the theory of stochastic differential equations, weak uniqueness
follows from trajectorial uniqueness, which is defined as follows. 
We
say that $(X,Y)$-trajectorial uniqueness on $\mathbb{E}_0$
holds for the equation
$\hat{F}_t(X,Y,Z)=0$  if 
equalities
$\hat{F}_t(\bx,\by,\bz)=0$ and
$\hat{F}_t(\bx,\by',\bz')=0$ for 
 $(\bx,\by,\bz)\in\mathbb{E}_0$ and 
$(\bx,\by',\bz')\in\mathbb{E}_0$ imply that 
 $\by=\by'$. 
In the next lemma, given $\bx$, we let $G(\bx)$ denote the set
of $\by$ such that $\hat{F}_t(\bx,\by,\bz)=0,\,t\in\R_+,$ for some
$\bz$, where $(\bx,\by,\bz)\in\mathbb{E}_0$. Note that
$G(\bx)=\emptyset$ if $\bx$ does not belong to the projection of
$\mathbb{E}_0$ on $\mathbb{E}^X$.
\begin{lemma}
  \label{le:weak}
If $(X,Y)$-trajectorial uniqueness  on $\mathbb{E}_0$ holds for the equation
$\hat{F}_t(X,Y,Z)=0$, then $\mathbf{\Pi}(Y=\by)=
\sup_{\bx\in G^{-1}(\by)}
\mathbf{\Pi}(X=\bx)$. In particular, $\mathbf{\Pi}(Y=\by)$ is
specified uniquely.
\end{lemma}
\begin{proof}
Let  $U$ be the subset of $\Upsilon$  such that
$F_t(X(\upsilon),Y(\upsilon),Z(\upsilon))=0$ and
$\bl(X(\upsilon),Y(\upsilon),Z(\upsilon)\br)\in\mathbb{E}_0$  for 
 $\upsilon\in U$. 
By hypotheses,  $\mathbf{\Pi}(\Upsilon\setminus U)=0$, so 
 $\mathbf{\Pi}(X=\bx,Y=\by)=\mathbf{\Pi}(\{X=\bx,Y=\by\}\cap
U)$. If $\upsilon\in U$ and $\mathbf{\Pi}(\upsilon)>0$, 
then the set $G(X(\upsilon))$ contains one
element.
 Hence, either
$\mathbf{\Pi}(X=\bx,Y=\by)=\mathbf{\Pi}(X=\bx)$ if
$\by\in G(\bx)$  or
$\mathbf{\Pi}(X=\bx,Y=\by)=0$ otherwise.
Consequently, $\mathbf{\Pi}(Y=\by)=\sup_{\bx\in\mathbb{E}^X}
\mathbf{\Pi}(X=\bx,Y=\by)=\sup_{\bx\in G^{-1}(\by)}
\mathbf{\Pi}(X=\bx)=\sup_{\bx\in G^{-1}(\by)}
\mathbf{\Pi}^X(\bx)$. 
\end{proof}
Thus, the proof of the theorem is completed by establishing 
trajectorial uniqueness.
In what follows, we  implement this programme.
\paragraph{Stochastic equations}
\label{sec:stochastic-equations}

By \eqref{eq:2}
\begin{multline}
  \label{eq:32}
  \overline{Q}^n_k(t)=\overline{Q}_k^n(0)+\sum_{m\in C_k}\Bl[
\int_0^t \ind\Bl(\frac{\overline{Q}^n_k(s-)}{w_{km}}
<\min_{\substack{l\in S_m:\\l\not=k}}\frac{\overline{Q}^n_{l}(s-)}{w_{lm}}\Br)
\,d\overline{A}^n_m(s)\\+
\int_0^t \ind\Bl(\frac{\overline{Q}^n_k(s-)}{w_{km}}
=\min_{\substack{l\in S_m:\\l\not=k}}\frac{\overline{Q}^n_l(s-)}{w_{lm}}\Br)
\alpha^n_{km}(ns)\,d\overline{A}^n_m(s)\Br]-
\int_0^t\ind\bl(\overline{Q}^n_k(s-)>0\br)\,d\overline{B}^n_k(s)\,.
\end{multline}
We introduce the following latent variables
\begin{align}
  \label{eq:5}
\overline{D}^n_{k}(t)&=  
\int_0^t\ind\bl(\overline{Q}^n_k(s-)>0\br)\,d\overline{B}^n_k(s),\\
    \overline{E}^n_{km}(t)&=
\int_0^t \ind\Bl(\frac{\overline{Q}^n_k(s-)}{w_{km}}<\min_{\substack{l\in S_m:\\l\not=k}}
\frac{\overline{Q}^n_l(s-)}{w_{lm}}\Br)\,d\overline{A}^n_m(s)\nonumber
\\&+  \label{eq:99}
\int_0^t
\ind\Bl(\frac{\overline{Q}^n_k(s-)}{w_{km}}
=\min_{\substack{l\in S_m:\\l\not=k}}\frac{\overline{Q}^n_l(s-)}{w_{lm}}
\Br)
\alpha^n_{km}(ns)\,d\overline{A}^n_m(s)\,.
\end{align}
They enable us to replace \eqref{eq:3} with a system of ``continuous''
 equations. 
 By \eqref{eq:32},   \eqref{eq:5}, and \eqref{eq:99}  for
 $t\in\R_+$ and $k=1,2,\ldots,K$
 \begin{subequations}
   \begin{align}
  \label{eq:3}
  \overline{Q}^n_k(t)&=\overline{Q}^n_k(0)
+\sum_{m\in C_k}\overline{E}^n_{k}(t)
-\overline{D}^n_k(t),\\\label{eq:3a}
  \int_0^t\overline{Q}^n_k(s-)\,d\overline{D}^n_k(s)&=
\int_0^t\overline{Q}^n_k(s-)\,d\overline{B}^n_k(s),\\
  \label{eq:109}
  \int_0^t \Bl(\frac{\overline{Q}^n_k(s-)}{w_{km}}-
\min_{l\in S_m}\frac{\overline{Q}^n_l(s-)}{w_{lm}}
\Br)
\,d\overline{E}^n_{km}(s)&=
  0,\\
  \label{eq:7}
\overline{A}^n_m(t)&=  \sum_{k\in S_m} \overline{E}^n_{km}(t).
\end{align}
 \end{subequations}
Let 
 $\overline{D}^n_k=\bl(\overline{D}^n_k(t),\,t\in\R_+\br)$
and $\overline{E}^n_{km}=\bl(\overline{E}^n_{km}(t),\,t\in\R_+\br)$,
where $k=1,2,\ldots,K$ and $m=1,2,\ldots,M$,
$\overline{D}^n=(\overline{D}^n_k,\,k=1,2,\ldots,K)$ and
$\overline{E}^n=\bl(\overline{E}^n_{km},\,k=1,2,\ldots,K,\,
m=1,2,\ldots,M\br)$. The latent
 variables $Z^n$ are defined by
$Z^n=\bl(\overline{E}^n,\,
\overline{D}^n\bl)$ and assume values in
$\mathbb{E}^Z=\D(\R_+,\R^{K\times M})\times
\D(\R_+,\R^K)$.
Below, generic elements of  spaces $\mathbb{E}^X$, $\mathbb{E}^Y$, and
 $\mathbb{E}^Z$  are denoted as follows: 
for space $\mathbb{E}^X$, as
  $\bx=(q_0,\ba,\bb)$, where $q_0=(q_{0,1},\ldots,q_{0,K})\in\R^K$,
$\ba=(\ba_m,\,m=1,2,\ldots,M)\in\D(\R_+,\R^M)$, 
 and
$\bb=(\bb_k,\,k=1,2,\ldots,K)\in\D(\R_+,\R^K)$,
 for space $\mathbb{E}^Y$, as
   $\by=\bq$, where
$\bq=(\bq_k,\,k=1,2,\ldots,K)\in\D(\R_+,\R^K)$,
and for space $\mathbb{E}^Z$, as
$\bz=(\be,\bd)$,
where $\be=(\be_k,\,k=1,2,\ldots,K)\in\D(\R_+,\R^K)$ and 
$\bd=(\bd_k,\,k=1,2,\ldots,K)\in\D(\R_+,\R^K)$.

\paragraph{Exponential tightness}
\label{sec:expon-tightn}

As pointed out above,
the hypotheses of the theorem 
 imply that  the $X^n$ LD converge in distribution at rate $n$ 
to an idempotent  variable $X=(Q_0,A,B)$ with idempotent distribution 
$\mathbf{\Pi}^X$. Since $A$ and $B$  have
 continuous paths $\mathbf{\Pi}$-a.e., 
the sequence $(A^n,B^n)$ is $\C$-exponentially tight.
Let us show that 
the sequence $(A^n,B^n,Y^n,Z^n)$  is $\C$-exponentially tight.
We repeat the argument of  the proof of Lemma~4.1 in Puhalskii \cite{Puh05}.
By \eqref{eq:32},  \eqref{eq:5}, and \eqref{eq:99} 
the increments of the $\overline{Q}^n_k(t)$, $\overline{E}^n_{km}(t)$
 and $\overline{D}^n_k(t)$ are majorised as follows by the increments of
 $\overline{A}^n_{m}(t)$ and $\overline{B}^n_k(t)$: for $s<t$
        \begin{align}
         \abs{\overline{Q}^n_k(t)-\overline{Q}^n_k(s)}&\le
     \sum_{m\in C_k}\abs{\overline{A}^n_{m}(t)-\overline{A}^n_{m}(s)}+
\nonumber     \abs{\overline{B}^n_k(t)-\overline{B}^n_k(s)},\\
     \abs{\overline{D}^n_k(t)-\overline{D}^n_k(s)}&\le
     \abs{\overline{B}^n_k(t)-\overline{B}^n_k(s)},     \label{eq:18}\\
     \abs{\overline{E}^n_{km}(t)-\overline{E}^n_{km}(s)}&\le
     \abs{\overline{A}^n_m(t)-\overline{A}^n_m(s)}.
\label{eq:181}
   \end{align}
 Since the sequences $\{\overline{A}^n_m,\,n\in\N\}$
 and
 $\{\overline{B}^n_k,\,n\in\N\}$  are $\C$-exponentially tight and the sequence
$\{\overline{Q}^n(0),\,n\in\N\}$ is exponentially tight, the above
inequalities imply by
\eqref{eq:c_exp_tight} that the sequences $\{\overline{Q}^n_k,\,n\in\N\}$,
$\{\overline{E}^n_{km},\,n\in\N\}$ and $\{\overline{D}^n_k,\,n\in\N\}$ 
 are $\C$-exponentially tight.  Therefore, 
 the sequences $Y^n$ and $Z^n$ are 
 $\C$-exponentially tight, so the sequence
 $\{(A^n,B^n,Y^n,Z^n),\,n\in\N\}$ is $\C$-exponentially tight.
As a consequence, 
the sequence $\{(X^n,Y^n,Z^n),\,n\in\N\}$ is exponentially tight.

\paragraph{Large deviation limit}
\label{sec:large-devi-limit-1}

Let the $(X^n,Y^n,Z^n)$ LD converge along a subsequence to an idempotent
process $(X,Y,Z)$ defined on an idempotent probability space
$(\Upsilon,\mathbf{\Pi})$,  where $X=( Q_0,  A, B)$, $Y= Q$, and
$Z=(E, D)$.  Our ultimate goal is to show that
$\mathbf{\Pi}^Q(\bq)=\mathbf{\Pi}(Q=\bq)$ is specified uniquely.
Note that $\mathbf{\Pi}$-a.e. $A(0)=0\in\R^M$, 
$B(0)=D(0)=0\in\R^K$, and $E(0)=0\in\R^{K\times M}$.
 Let us show
that $\mathbf{\Pi}$-a.e. the idempotent processes $D_k$ are absolutely
continuous, nonnegative and
nondecreasing. Since
$\mathbf{P}(\overline{D}^n_k(t)-\overline{D}^n_k(s)<0)=0$ for $t\ge
s$ and $\liminf_{n\to\infty}
\mathbf{P}(\overline{D}^n_k(t)-\overline{D}^n_k(s)<0)^{1/n}\ge
\mathbf{\Pi}( D_k(t)- D_k(s)<0)$, we obtain
that $\mathbf{\Pi}( D_k \text{ is not nondecreasing})=\sup_{s<t}
\mathbf{\Pi}( D_k(t)- D_k(s)<0)=0$, so
$ D_k$ is nondecreasing $\mathbf{\Pi}$-a.e. It is therefore nonnegative.
Next, by  \eqref{eq:18}
$  \mathbf{P}\bl(\abs{ \overline{D}^n_k(t)-\overline{D}^n_k(s)}>
     \abs{\overline{B}^n_k(t)-\overline{B}^n_k(s)}\br)=0$,
so similarly to the above
 $\mathbf{\Pi}(\abs{ D_k(t)- D_k(s)}>
     \abs{ B_k(t)- B_k(s)})=0$.
Since the function $ B_k$ is absolutely continuous
$\mathbf{\Pi}$-a.e., 
it follows that so is $ D_k$ and
$   \dot{D}_k(t)\le \dot{B}_k(t).$
A similar argument applied to the $ E_{km}$ and using
\eqref{eq:181} shows that the $ E_{km}$ are
$\mathbf{\Pi}$-a.e. absolutely continuous and nondecreasing.

Letting $n\to\infty$ in \eqref{eq:3}, \eqref{eq:3a}, \eqref{eq:109}
and \eqref{eq:7} we obtain by the continuous mapping principle 
that $\mathbf{\Pi}$-a.e.
 \begin{subequations}
   \begin{align*}
   Q_k(t)&= Q_{0,k}+\sum_{m\in C_k} E_{km}(t)
- D_k(t),\\
  \int_0^t Q_k(s)\,d D_k(s)&=
\int_0^t Q_k(s)\,d B_k(s),\\
  \int_0^t \Bl(
\frac{ Q_k(s)}{w_{km}}-\min_{l\in S_m}\frac{ Q_l(s)}{w_{lm}}\Br)
\,d E_{km}(s)&=0,\\
 A_m(t)&=  \sum_{k\in S_m}  E_{km}(t).
\end{align*}
 \end{subequations}
\paragraph{Trajectorial uniqueness}
\label{sec:traj-uniq}

We take as  $\mathbb{E}_0$ 
 the subset of $\R^K\times\D(\R_+,\R^{3K+M+KM+1})$ of
elements $(q_0,\bq,\ba,\bb,\be,\bd)$ such that $\bq(0)=q_0$,
the functions $\bq,\ba,\bb,\bd$ and $\be$ are componentwise nonnegative and
absolutely continuous, the functions $\ba,\bb,\bd$ and $\be$
 are nondecreasing and start at $0$, 
 and the following relations hold
for $k=1,2,\ldots,K,\,m=1,2,\ldots,M,\,t\in\R_+$
\begin{subequations}
  \begin{align}
  \label{eq:1}
  \dot{\bq}_k(t)=\sum_{m\in C_k}\dot{\be}_{km}(t)-
\dot{\bd}_k(t),\;\bq_k(0)=q_{0,k},\\
      \label{eq:159}
\bq_k(t)\bl(\dot{\bd}_k(t)-\dot{\bb}_k(t)\br)=0,\\
\label{eq:13}
\Bl(\frac{\bq_k(t)}{w_{km}}-\min_{l\in S_m}\frac{\bq_{l}(t)}{w_{lm}}
\Br)\dot{\be}_{km}(t)=0,\\
  \label{eq:8}
\dot{\ba}_m(t)=\sum_{k\in S_m}\dot{\be}_{km}(t),
\;\dot{\bd}_k(t)\le\dot{\bb}_k(t).
    \end{align}
\end{subequations}
We have proved that  
$\mathbf{\Pi}\bl((Q_0,A,B,Q,E,D)\not\in\mathbb{E}_0\br)=0$. 
Hence, $\mathbf{\Pi}(Q=\bq)=0$ if $\bq$ either is not absolutely
continuous, or componentwise nonnegative, 
 or it is not obtained as a solution
of the latter system of equations.
We next prove that this system of equations 
uniquely specifies $\bq$ given $q_0$, $\ba$ and $\bb$. 

 Let $\bq'$ be another solution, i.e.,
\begin{align*}
    \dot{\bq}'_k(t)=\sum_{m\in C_k}\dot{\be}'_{km}(t)-
\dot{\bd}'_k(t),\;\bq'_k(0)=q_{0,k},\\
\bq'_k(t)\bl(\dot{\bd}'_k(t)-\dot{\bb}_k(t)\br)=0,\\
\Bl(\frac{\bq'_k(t)}{w_{km}}-\min_{l\in S_m}\frac{\bq'_{l}(t)}{w_{lm}}
\Br)\dot{\be}'_{km}(t)=0,\\
\dot{\ba}_m(t)=\sum_{k\in S_m}\dot{\be}'_{km}(t),\;
\dot{\bd}'_k(t)\le\dot{\bb}_k(t).
\end{align*}

It suffices to prove that 
\begin{equation}
  \label{eq:9}
  \frac{d}{dt}\sum_{k=1}^K\abs{\bq_k(t)-\bq'_k(t)}\le0.
\end{equation}
  We have 
  \begin{multline}
    \label{eq:19}
  \frac{d}{dt}\sum_{k=1}^K\abs{\bq_k(t)-\bq'_k(t)}
=\sum_{k=1}^K\Bl(\bl(\dot{\bq}_k(t)-\dot{\bq}'_k(t)\br)
\ind(\bq_k(t)>\bq'_k(t))\\+
\bl(\dot{\bq}'_k(t)-\dot{\bq}_k(t)\br)
\ind(\bq_k(t)<\bq'_k(t))\Br)\,.   
  \end{multline}
By \eqref{eq:1}
\begin{equation}
  \label{eq:15}
  \dot{\bq}_k(t)-\dot{\bq}'_k(t)=
\sum_{m\in C_k}\bl(\dot{\be}_{km}(t)-\dot{\be}'_{km}(t)\br)+
\bl(\dot{\bd}'_k(t)-\dot{\bd}_k(t)\br).
\end{equation}
If $\bq_k(t)>\bq'_k(t)$, then $\bq_k(t)>0$, so by
\eqref{eq:159} $\dot{\bd}_k(t)=\dot{\bb}_k(t)$, which implies by 
the inequality $\dot{\bd}'_k(t)\le\dot{\bb}_k(t)$  and \eqref{eq:15} that
$\dot{\bq}_k(t)-\dot{\bq}'_k(t)\le
\sum_{m\in C_k}\bl(\dot{\be}_{km}(t)-\dot{\be}'_{km}(t)\br)$ a.e
on the set $\{\bq_k(t)>\bq'_k(t)\}$. Similarly,
$\dot{\bq}'_k(t)-\dot{\bq}_k(t)\le
\sum_{m\in C_k}\bl(\dot{\be}'_{km}(t)-\dot{\be}_{km}(t)\br)$
 a.e on the set $\{\bq'_k(t)>\bq_k(t)\}$.
Hence, by \eqref{eq:19} 
\begin{multline*}
  \frac{d}{dt}\sum_{k=1}^K\abs{\bq_k(t)-\bq'_k(t)}
\le\sum_{m=1}^M \sum_{k\in
  S_m}\Bl(\bl(\dot{\be}_{km}(t)-\dot{\be}'_{km}(t)\br)
\ind(\bq_k(t)>\bq'_k(t))\\+\bl(\dot{\be}'_{km}(t)-\dot{\be}_{km}(t)\br)
\ind(\bq_k(t)<\bq'_k(t))\Br)\,.
\end{multline*}
We prove that on the  right-hand side each sum over $S_m$ is
nonpositive.

If $  \min_{l\in S_m}\bq_l(t)/w_{lm}= \min_{l\in
  S_m}\bq'_l(t)/w_{lm}$, then on the set where $\bq_k(t)>\bq'_k(t)$ we
have that   $\bq_k(t)/w_{km}>\bq'_k(t)/w_{km}\ge
\min_{l\in S_m}\bq_l(t)/w_{lm}$,  so by \eqref{eq:13} 
$\dot{\be}_{km}(t)=0$. Similarly, $\dot{\be}'_{km}(t)=0$ a.e. on the
  set where $\bq'_k(t)>\bq_k(t)$. The required property follows.

If 
$  \min_{l\in S_m}\bq_l(t)/w_{lm}< \min_{l\in
  S_m}\bq'_l(t)/w_{lm}$, then, analogously to the preceding argument, 
 $\dot{\be}_{km}(t)=0$ a.e.
 on the set where $\bq_k(t)\ge\bq'_k(t)$, and
the sum in question is not greater than
$\sum_{k\in S_m}\bl(\dot{\be}'_{km}(t)-\dot{\be}_{km}(t)\br)=
\dot{\ba}_m(t)-\dot{\ba}_m(t)=0$.
Similarly, if $  \min_{l\in S_m}\bq_l(t)/w_{lm}> \min_{l\in
  S_m}\bq'_l(t)/w_{lm}$, then  $\dot{\be}'_{km}(t)=0$
on the set where $\bq_k(t)\le\bq'_k(t)$ so that
this sum  is not greater than
$\sum_{k\in S_m}\bl(\dot{\be}_{km}(t)-\dot{\be}'_{km}(t)\br)=0$.
Inequality \eqref{eq:9} has been proved.
\paragraph{Evaluating the limit idempotent distribution}
\label{sec:eval-acti-funct}

By
Lemma~\ref{le:weak}
\begin{equation}
  \label{eq:4}
\mathbf{\Pi}(Q=\bq)=\sup_{(q_0,\ba,\bb)\in\Delta(\bq)}
\mathbf{\Pi}(Q_0=q_0,A=\ba,B=\bb),
\end{equation}
where $\Delta(\bq)$ is the set of  $(q_0,\ba,\bb)$ such that
\eqref{eq:1}, \eqref{eq:159}, \eqref{eq:13} and \eqref{eq:8}
 are satisfied.

We conclude the proof by evaluating
  the right-hand side of \eqref{eq:4}. Let 
$\mathbf{I}(\bq)=-\log \mathbf{\Pi}(Q=\bq)$. Then
\begin{equation}
  \label{eq:23}
  \mathbf{I}(\bq)=\inf_{(q_0,\ba,\bb)\in\Delta(\bq)}
\bl(\mathbf{I}^{Q_0}(q_0)+\mathbf{I}^{A,B}(\ba,\bb)\br)=
\mathbf{I}^{Q_0}(\bq(0))+\inf_{(\ba,\bb)\in\tilde{\Delta}(\bq)}
\int_0^\infty\psi(\dot{\ba}(t),\dot{\bb}(t))
\,dt,
\end{equation}
where $\tilde{\Delta}(\bq)$ is the set of  $(\ba,\bb)$ such that
\eqref{eq:1}, \eqref{eq:159}, \eqref{eq:13} and \eqref{eq:8}
 are satisfied with
$\bq(0)$ as $q_0$. The definition of $L(x,y)$ implies that 
$\psi(\dot{\ba}(t),\dot{\bb}(t))\ge
L(\bq(t),\dot{\bq}(t))$ a.e. for $(\ba,\bb)\in\tilde{\Delta}(\bq)$.
Hence,
$  \mathbf{I}(\bq)\ge \mathbf{I}^Q(\bq).$

In order to prove the reverse inequality, we assume that
$\mathbf{I}^Q(\bq)<\infty$, so $L(\bq(t),\dot{\bq}(t))<\infty$ a.e. 
Since  the infimum in the definition \eqref{eq:12} of $L(x,y)$ is
attained if finite, the sets  
    $\Gamma(t),  t\in \R_+$, consisting of 
 $(a,b)\in N(\bq(t),\dot{\bq}(t))$
such that $\psi(a,b)=L(\bq(t),\dot{\bq}(t))$ are nonempty for
almost all $t$.
Also, the graph $\{(t,(a,b))\in \R_+\times
\R^{M+K}:\,(a,b)\in\Gamma(t)\}$ is measurable with
respect to the product of the Lebesgue $\sigma$-algebra on $\R_+$ and
Borel $\sigma$-algebra on $\R^{M+K}$. Therefore, by a measurable
selection theorem (see, e.g., Clark \cite{Cla83}) there exist Lebesgue
measurable functions $\tilde{\ba}$ and $\tilde{\bb}$
such that $(\tilde{\ba}(t),\tilde{\bb}(t))\in \Gamma(t)$
 a.e. Letting $\ba(t)=\int_0^t \tilde{\ba}(s)\,ds$  and 
$\bb(t)=\int_0^t \tilde{\bb}(s)\,ds$, we obtain that
$(\ba,\bb)\in\tilde{\Delta}(\bq)$ and 
$\int_0^\infty\psi(\dot{\ba}(t),\dot{\bb}(t))
\,dt=\int_0^\infty L(\bq(t),\dot{\bq}(t))\,dt.$
Thus the infimum on the right of \eqref{eq:23} is equal to
$\int_0^\infty L(\bq(t),\dot{\bq}(t))\,dt$, hence,
 $\mathbf{I}(\bq)=\mathbf{I}^Q(\bq)$ and
$\mathbf{\Pi}^Q(\bq)=\exp\bl(-\mathbf{I}^Q(\bq)\br)$.
\paragraph{Acknowledgements}
\label{sec:acknowledgements}

The first author is thankful to Yuri Suhov for drawing his attention to the
join-the-shortest-queue setting.
\appendix

\section{Review of idempotent probability and large deviation convergence}
\label{sec:backgr-idemp-prob}

Let $\Upsilon$ be a set. A  function $\mathbf{\Pi}$ from the power
set of $\Upsilon$ to $[0,1]$ is called an idempotent probability
 if $\mathbf{\Pi}(\Gamma)=\sup_{\upsilon\in\Gamma}
\mathbf{\Pi}(\{\upsilon\}),\,\Gamma\subset \Upsilon$ and
$\mathbf{\Pi}(\Upsilon)=1$. The pair $(\Upsilon,\mathbf{\Pi})$ is
called an idempotent probability space.
For economy of notation, 
we denote $\mathbf{\Pi}(\upsilon)=\mathbf{\Pi}(\{\upsilon\})$.
A property $\mathcal{P}(\upsilon),\,\upsilon\in\Upsilon,$ 
pertaining to the elements of $\Upsilon$ is said to hold
$\mathbf{\Pi}$-a.e. if $\mathbf{\Pi}(\mathcal{P}(\upsilon)
\text{ does not hold})=0$.
A function $f$ from a set
$\Upsilon$ equipped with 
idempotent probability $\mathbf{\Pi}$ to a set $\Upsilon'$ 
is called an idempotent variable. 
 The idempotent distribution of an idempotent variable $f$ is defined
 as the set function $\mathbf{\Pi}\circ
f^{-1}(\Gamma)=\mathbf{\Pi}(f\in\Gamma),\,\Gamma\subset \Upsilon'$.
If $f$ is the canonical idempotent variable that is defined by
$f(\upsilon)=\upsilon$, then it has $\mathbf{\Pi}$ 
as the idempotent distribution.
$\Upsilon'$-valued idempotent 
variables $f$ and $f'$ are said to be independent if 
$\mathbf{\Pi}(f=\upsilon',\,f'=\upsilon'')=\mathbf{\Pi}(f=\upsilon')
\mathbf{\Pi}(f'=\upsilon'')$ for all $\upsilon',\upsilon''\in
\Upsilon'$. 
Independence of finite collections of idempotent variables is defined
similarly. 
A collection $(X_t,\,t\in\R_+)$ 
of $\R^\ell$-valued idempotent variables on $\Upsilon$, where $\ell\in\N$,
is called an idempotent process.
The functions $(X_t(\upsilon),\,t\in\R_+)$ 
for various $\upsilon\in\Upsilon$ are
called trajectories (or paths) of $X$.
Idempotent processes are said to be independent if they are
independent as idempotent variables with values in the associated
function space.

 If $\Upsilon$ is, in addition, a metric space and 
 the sets $\{\upsilon\in\Upsilon:\,\mathbf{\Pi}(\upsilon)\ge \alpha\}$ 
are compact for all
$\alpha\in(0,1]$, then $\mathbf{\Pi}$ is called a deviability.  
Obviously, $\mathbf{\Pi}$ is a deviability if and only if
$\mathbf{I}(\upsilon)=-\log\mathbf{\Pi}(\{\upsilon\})$ is an action
 functional.
 If $f$ is a continuous mapping from
$\Upsilon$ to another metric space $\Upsilon'$, then 
 $\mathbf{\Pi}\circ f^{-1}$ is a deviability on $\Upsilon'$.
As a matter of fact, for the latter property to hold, 
one can only require that $f$ be continuous on
the sets $\{\upsilon\in\Upsilon:\,\mathbf{\Pi}(\upsilon)\ge \alpha\}$ for
 $\alpha\in(0,1]$.
In general,  $f$ is said to be  a Luzin idempotent variable 
if  $\mathbf{\Pi}\circ f^{-1}$ is a
deviability on $\Upsilon'$. 

Let  $\{\mathbf{P}_n,\,n\in\N\}$ be a  sequence of probability measures on
a metric space $\mathbb{E}$ endowed with  Borel $\sigma$-algebra and
let $\mathbf{\Pi}$ be a deviability on $\mathbb{E}$.
 Let $m_n\to\infty$ as $n\to\infty$. 
The sequence $\{\mathbf{P}_n,\,n\in\N\}$
 is said to large deviation converge (LD converge) at rate $m_n$
 to  $\mathbf{\Pi}$ as $n\to\infty$ if $\lim_{n\to\infty}\Bl(\int_\mathbb{E}
 f(z)^{m_n}\,\mathbf{P}_n(dz)\Br)^{1/m_n}= 
\sup_{z\in \mathbb{E}}f(z)\mathbf{\Pi}(z)$
 for every bounded continuous $\R_+$-valued function $f$ on $\mathbb{E}$.
This definition is equivalent to requiring that the inequalities 
$\limsup_{n\to\infty}\mathbf{P}_n(F)^{1/m_n}\le \mathbf{\Pi}(F)$ and 
$\liminf_{n\to\infty}\mathbf{P}_n(G)^{1/m_n}\ge \mathbf{\Pi}(G)$ hold for all
closed sets $F$ and all open sets $G$, respectively.
Therefore, the sequence $\{\mathbf{P}_n,\,n\in\N\}$ LD 
converges at rate $m_n$ to
 $\mathbf{\Pi}$ if and only if it obeys the LDP for scale $m_n$ 
with action functional 
 $\mathbf{I}(z)=-\log\mathbf{\Pi}(z)$. We favour  the term ``LD
 convergence'' over ``the LDP'' as being more natural for our
 approach.   The deviability $\mathbf{\Pi}$
 is said to be an LD limit point of the $\mathbf{P}_n$ for rate $n$ if
 each subsequence  $\{\mathbf{P}_{n_k},\,k\in\N\}$ 
of $\{\mathbf{P}_n,\,n\in\N\}$
 contains a further subsequence $\{\mathbf{P}_{n_{k_l}},\,l\in\N\}$ that LD
 converges to $\mathbf{\Pi}$ at rate $n_{k_l}$ as $l\to\infty$.
The sequence $\{\mathbf{P}_n,\,n\in\N\}$  is said to be exponentially tight on
order $m_n$ if for arbitrary $\epsilon>0$ there exists a compact subset
$E$ of $\mathbb{E}$ such that
$  \limsup_{n\to\infty}\mathbf{P}_n(\mathbb{E}\setminus E)^{1/m_n}<\epsilon.$
An exponentially tight sequence possesses LD limit
points. Thus, one can prove LD convergence of the $\mathbf{P}_n$ 
by proving that exponential
tightness holds and that there is a unique LD limit point.
We will say that a sequence of random variables with values in a
metric space is exponentially tight
if so is the sequence of their laws.

LD convergence of probability measures 
can be also expressed as LD convergence in distribution of the
associated random variables to idempotent variables. 
In the setting of stochastic processes,
this point of view enables one to
consider the LD limit as a dynamical system rather than as ``a mass
function'' on the space of trajectories.  We say that 
 a sequence
$\{X_n,\,n\in\N\}$  of random variables defined on probability spaces
$(\Omega_n,\mathcal{F}_n,\mathbf{P}_n)$ and assuming  values in $\mathbb{E}$ 
LD converges in distribution at rate $m_n$ as $n\to\infty$
to a Luzin idempotent variable $X$ defined on an 
idempotent probability space $(\Upsilon,\mathbf{\Pi})$ and assuming
values in $\mathbb{E}$ if the sequence of the probability laws of the $X_n$ LD
converges to the idempotent distribution of $X$ at rate $m_n$. 
Conversely, LD convergence of a sequence $\{\mathbf{P}_n,\,n\in\N\}$
of probability measures on $\mathbb{E}$ to a deviability
$\mathbf{\Pi}$ on $\mathbb{E}$ can expressed as LD convergence in
distribution if one considers the canonical setting.
 The continuous mapping principle, known as the
contraction principle for the LDP, states that if the
$X_n$ LD converge in distribution to $X$ and $f$ is a continuous
function from $\mathbb{E}$ to another metric space, then the $f(X_n)$ LD
converge in distribution to $f(X)$. We will use the extension, also
referred to as the continuous mapping principle, where the function
$f$ is allowed to be a measurable function that is continuous only
a.e. with respect to the idempotent distribution of $X$. For a
detailed discussion, see Garcia \cite{Gar04}.
  The definition of a  limit point for LD
convergence in distribution is similar to that for LD convergence of
probability measures.

If the $X_n$ are stochastic processes with trajectories in a Skorohod
 space $\D(\R_+,\R^\ell)$, then the sequence 
 $X_n$ is said to be $\C$-exponentially tight on order $m_n$
 if it is exponentially  tight on order $m_n$ and each LD limit
 point $\mathbf{\Pi}$ 
of the distributions of the $X_n$ is an idempotent distribution
 of a continuous-path idempotent process in the sense that 
$\mathbf{\Pi}\bl(\D(\R_+,\R^\ell)\setminus \C(\R_+,\R^\ell)\br)=0$.
The sequence 
$X_n$ is $\C$-exponentially tight on order $m_n$ if and
only if
\begin{equation}\label{eq:c_exp_tight}
  \begin{split}
    \lim_{L\to\infty}\limsup_{n\to\infty}
\mathbf{P}\bl(\abs{X_n(0)}>L\br)^{1/m_n}&=0,\\
  \lim_{\delta\to0}\limsup_{n\to\infty}
\mathbf{P}\bl(\sup_{\substack{s,t\in[0,T]:\\\abs{s-t}\le\delta}}
\abs{X_n(t)-X_n(s)}>\epsilon\br)^{1/m_n}&=0,\;T\in\R_+,\,\epsilon>0.
  \end{split}
\end{equation}
If the sequence 
 $X_n$ is  $\C$-exponentially tight, then a limit point of the $X_n$
 for LD convergence in distribution may be  considered as an
 idempotent process with trajectories in $\C(\R_+,\R^\ell)$.
Converesly,  if the $X_n$ LD converge in distribution to a
continuous-path idempotent process, then the sequence $X_n$ is
$\C$-exponentially tight.

We also note that, as a consequence of the continuous mapping
 principle, if a sequence of stochastic processes $(X_n,Y_n)$ assuming
 values in $\R^{\ell_1}$ and $\R^{\ell_2}$, respectively,
 LD converges in distribution in   $\D(\R_+,\R^{\ell_1})\times
 \D(\R_+,\R^{\ell_2})$ to a continuous-path
 idempotent process, then the LD convergence also holds in 
$\D(\R_+,\R^{\ell_1}\times\R^{\ell_2})$.

\def\cprime{$'$} \def\cprime{$'$}

\end{document}